\documentclass[12pt,a4paper]{article}
\usepackage{amsmath}

\usepackage{amsfonts}
\usepackage{amssymb}
\usepackage{amsmath}
\usepackage{amsthm}
\usepackage{amscd}
\usepackage{graphicx}
\usepackage{color}
\newcommand{\be}{\begin{equation}}
\newcommand{\ee}{\end{equation}}
\newcommand{\ef}[1]{\, #1}

\newcommand\rd{{\mathrm d}}

\newcommand\fois{\mathord{\cdot}}
\newcommand\dd{{\text{\textup{d}}}}

\newcommand\caJ{{\mathcal J}}

\newcommand\caP{{\mathcal P}}

\newcommand\gC{{\mathbb C}}

\newcommand\gN{{\mathbb N}}

\newtheorem{Theorem}{Theorem}[section]
\newtheorem{theorem}[Theorem]{Theorem}

\newtheorem{conj}[Theorem]{Conjecture}

\newtheorem{definition}[Theorem]{Definition}

\newcommand{\beqa}{\begin{eqnarray}}
\newcommand{\eeqa}{\end{eqnarray}}

\usepackage[latin1]{inputenc}
\usepackage{subfigure}

\author{A. Tanasa\thanks{adrian.tanasa@labri.fr, Partially supported by the PN  09370102 grant}}
\title
{Combinatorial quantum~field~theory and the Jacobian conjecture}
\begin{document}
\maketitle
\begin{abstract}
In this short review we first recall combinatorial 
or ($0-$dimensional) quantum field theory (QFT). We then give the main idea of a standard QFT method, called
the intermediate field method, and we 
review how to apply this method to a combinatorial QFT reformulation of the celebrated Jacobian Conjecture
on the invertibility of polynomial systems.
This approach establishes a related theorem concerning partial elimination of variables that implies a reduction of the generic case to the quadratic one. Note that this does not imply solving the Jacobian Conjecture, because one needs to introduce a supplementary parameter for the dimension of a certain linear subspace where the system holds. 
\end{abstract}




\section{Introduction}
\label{sec:in}

In 1939, Keller formulated in \cite{Keller} the \textbf{Jacobian Conjecture},
as a remarkably simple and natural conjecture on the 
invertibility of polynomial systems.
The conjecture states that
 a polynomial function $F : {\mathbb{C}}^n\to   {\mathbb{C}}^n$ 
 is invertible, and its inverse is polynomial, if and
only if the determinant of the Jacobian matrix 
$J_F(z) = \big( \rd F_j(z)/\rd z_i \big)_{1 \leq i,j \leq n}$ is a
non-zero constant.\footnote{In \cite{Keller}, Keller formulated this conjecture for $n=2$ and polynomials with integral coefficients, but this conjecture was shortly after generalized to the form above.}

In order to illustrate this celebrated conjecture, let us give a simple example. Let
$n=2$, $F(z_1, z_2)=(z_1+z_2^3, z_2)$. 
The 
Jacobian writes 
 \beqa
{\mathrm{det}} 
\begin{pmatrix}
\frac{d F_1}{dz_1} & \frac{d F_1}{dz_2}\\
\frac{d F_2}{dz_1} & 
\frac{d F_2}{dz_2}
\end{pmatrix}
= det  
\begin{pmatrix}
1 & 3z_2^2 \\
0 & 1
\end{pmatrix}
 = 1
 \eeqa
One easily finds: $F^{-1}(z_1,z_2)=(z_1-z_2^3,z_2)$. 

\medskip

Despite several efforts (of mathematicians such as Gr\"obner, Zariski, Oda, to namely only a few), and various promising
partial results, the Jacobian Conjecture remains unsolved.  An introduction to the problem,
the context, and the state of the art up to 1982, can be found in the
seminal paper
\cite{BassCW}.

The function $F: \gC^n \to \gC^n$ is said to be a 
\emph{polynomial system} if all the coordinate functions $F_j$'s are polynomials. Let us
call $\caP_n$ the set of such functions. For a function $F$, define
\be
\label{def.jac}
J_F(z) = 
\left(
\frac{\rd}{\rd z_i} F_j(z) \right)_{1 \leq i,j \leq n},
\ee
the corresponding Jacobian matrix
and the two subspaces of~$\caP_n$:

\begin{definition}
\begin{align}
\caJ_n^{\rm lin}
&:= \{F\in \caP_n \,|\, \det J_F(z) = c \in \gC^{\times}\},\\
\caJ_n
&:=\{F\in \caP_n\, |\, F {\text{ is invertible}}\}.
\end{align}
\end{definition}
It will often be sufficient to analyse the
subset of $\caJ^{\rm lin}$ such that $\det J_F(z) = 1$.
The Jacobian Conjecture rewrites as
\begin{conj}[Jacobian Conjecture, \cite{Keller}]
\label{conj.jc}
\be
\caJ_n^{\rm lin} = \caJ_n
\qquad
\forall n
\ef.
\ee
\end{conj}
Define the \emph{total degree} of $F$, $\deg(F)$, as $\max_j
\deg(F_j(z))$, and introduce the subspaces
\begin{align}
\caP_{n,d} &= \{ F \in \caP_n \;|\; \deg(F) \leq d \}
\ef;
\end{align}
and similarly for $\caJ$ and $\caJ^{\rm lin}$. 

Let us now mention two positive
results on the conjecture
\begin{enumerate}
\item 
a theorem for the quadratic case
($d=2$), established first in \cite{Wang}, and then
in \cite{Oda} (see also \cite[Lemma 3.5]{Wright2} and
\cite[Thm.\ 2.4, pag.\ 298]{BassCW}).
\begin{theorem}[\cite{Wang}]
\label{thm.red2}
\be
\caJ_{n,2}^{\rm lin} = \caJ_{n,2}
\qquad
\forall n
\ef.
\ee
\end{theorem}
\item a reduction theorem, from the general case to the cubic case, established by
Bass {\it et. al.} \cite[Sec.~II]{BassCW}).
\begin{theorem}[\cite{BassCW}]
\label{thm.red3}
\be
\caJ_{n,3}^{\rm lin} = \caJ_{n,3}
\quad
\forall n
\qquad
\Longrightarrow
\qquad
\caJ_{n}^{\rm lin} = \caJ_{n}
\quad
\forall n
\ef.
\ee
\end{theorem}
\end{enumerate}

%



The result 
of \cite{GST} that we will review here from a quantum field theoretical (QFT) perspective 
is a reduction theorem, somewhat analogous to the one
above, that tries to "fill the gap" between the cases of degree two
and three.  However,  we need to introduce an adaptation of the
statement of the Jacobian conjecture with one more parameter.
\begin{definition}
\label{def-jjlin}
For $n'\leq n$ and $F\in\caP_{n,d}$, we write $z=(z_1,z_2)$ and
$F=(F_1,F_2)$ to distinguish components in the two subspaces
$\gC^{n'}\times\gC^{n-n'} \equiv \gC^{n}$.
We set $R(z_2;z_1)$ as a synonymous for $F_2(z_1,z_2)$,
emphasizing that, in $R$, we consider $z_2$ as the variables in a
polynomial system, and $z_1$ as parameters. Then, we define the
subspaces of $\caP_{n,d}$
\begin{align*}
\caJ_{n,d;n'} & :=\{F\in\caP_{n,d}\,|\,
F^{-1}\text{~is defined on~}
\gC^{n'}\times\{0\}
\}
\\
\caJ^{\rm lin}_{n,d;n'} & :=\{F\in\caP_{n,d}\,|\,
R\in\caJ_{n-n',d}\text{ and }(\det J_F)(z_1,R^{-1}(0,z_1))=c \in \mathbb{C}^{\times},\;\forall z_1\in\gC^{n'}\}
\end{align*}
\end{definition}



Let us now state the reduction theorem we will review in this paper:
\begin{Theorem} [\cite{GST}]
\label{lem-main}
For $n\in\gN$ and $d\geq 3$, there exists an injective map 
$\Phi: \caP_{n,d}\mapsto \caP_{n(n+1),d-1}$ 
satisfying 
\begin{align}
\Phi(\caJ^{\rm lin}_{n,d})
&\equiv
\caJ^{\rm lin}_{n(n+1), d-1;n}\cap\Phi(\caP_{n,d})
\ef;
&
\Phi(\caJ_{n,d})
&\equiv
\caJ_{n(n+1), d-1;n}\cap\Phi(\caP_{n,d})
\ef.
\end{align}
\end{Theorem}
Combining Theorem \ref{thm.red3} and the theorem above, the full
Jacobian Conjecture reduces to the question whether
\beqa
\caJ^{\rm lin}_{n(n+1),2;n}=\caJ_{n(n+1),2;n}.
\eeqa







\section{Combinatorial (or $0-$dimensional) QFT and the intermediate field method}
\label{sec:OdimQFT}

In this section we introduce combinatorial (or $0-$dimensional) QFT and we then give the general idea of the so-called intermediate field method.

\subsection{Combinatorial (or $0-$dimensional) QFT}
\label{subsec:OdimQFT}

Usually in QFT, the scalar {field $ \varphi$} is function of space (${\mathbb{R}}^D$), see for example the review papers \cite{malek-SLC} or \cite{adrian-SLC}, and $D=4$.

However, one can consider the case 
$D=0$ .
In this case the scalar field $ \varphi_i$ is not a function of space
(since there is no space anymore),
and 
$\varphi$ is simply a (real or complex) variable.
Some authors refer to this simpler formulation of QFT as {\bf combinatorial QFT}, and we will do the same here.

Thus, in the real case, the 
{\it partition function} (which, from a combinatorial point of view, is a generating function) 
of the so-called $\varphi^4$ model 
writes as the integral
\beqa
\label{act:phi4}
Z= \int_{{\mathbb R}} d\varphi \, e^{-\frac 12 \varphi^2+\frac{\lambda}{4!} \varphi^4},
\eeqa
where the constant 
{$\lambda$ is the QFT coupling constant}.
The term $\varphi^4$ is called an interaction term of degree four.

It is worth emphasizing here that in combinatorial QFT,  functional integrals 
(which are particularly involved to rigorously define for $D\ge 1$)
become usual (real or complex)  integrals. Combinatorial QFT is thus much easier to manipulate from a mathematical point of view.

However, combinatorial QFT still presents a certain interest for the mathematical physicist
because it can be seen as some kind of "laboratory" to test the usual QFT mathematical tools.
Thus, one needs to evaluate integrals of the same type as in usual QFT, namely type integrals 
\beqa
\frac{\lambda^n}{n} \int d\varphi\, e^{-\varphi^2/2} \left( \frac{\varphi^4}{4!}\right)^n.
\eeqa
coming from the so-called QFT perturbative expansion
(which comes to Taylor expand the exponential in \eqref{act:phi4} and then dealing with each term one by one, instead of dealing with the integral as a whole).
In order to evaluate this type of integrals, 
one can use standard QFT techniques. Namely, one can define
\beqa
Z_0(J)=\int_{{\mathbb R}} d\varphi\,  e^{-\varphi^2/2 + J \phi}
\eeqa
where the $J$ is the QFT source.

Using this QFT source technique, the $(2k)-$point correlation functions can be computed in the following way:
\beqa
\int_{{\mathbb R}} d\phi\, e^{-\phi^2/2} {\varphi^{2k}}=
\frac{\partial^{2k}}{\partial J^{2k}}
\int_{{\mathbb R}} d\varphi\,  e^{-\varphi^2/2 + J \varphi}|_{J=0}=
\frac{\partial^{2k}}{\partial J^{2k}} e^{J^2/2}|_{J=0}.
\eeqa










\subsection{The intermediate field method}
\label{sec:inter}

We give now the general idea of the so-called  intermediate field method. Note that we present this method in the case of combinatorial QFT, but this idea can generalize to arbitrary $D$.

Thus, the intermediate field method consists of  introducing a new field, $\sigma$, used to rewrite the interaction
in a way that allows for 
the degree of the interaction to be reduced.


In order to illustrate this, let us take the example of the 
 $\varphi^6$ model, where $\varphi$ is, as above, a real $0-$dimensional field (thus, a real variable).


The partition function of the $\varphi^6$ combinatorial QFT writes model
\beqa
\label{act:phi6}
Z(\lambda)=\int _{{\mathbb R}}\frac{d\varphi}{\sqrt{2\pi}}e^{-\frac 12 \varphi^2}e^{-\lambda\varphi^6}
\eeqa

The intermediate field model consists of rewriting this partition function with the help of a supplementary integral on the intermediate field $\sigma$ in the following way:
\beqa
Z(\lambda)=\int_{{\mathbb R}} \frac{d\varphi}{\sqrt{2\pi}}e^{-\frac 12 \varphi^2}
\int_{{\mathbb R}} \frac{d\sigma}{\sqrt{2\pi}}e^{-\frac 12 \sigma^2}
e^{\imath \sqrt{2\lambda}\varphi^3\sigma}.
\eeqa
Note that this allowed to replace the interacting term $\varphi^6$ of 
\eqref{act:phi6} by a lower degree interacting term $\varphi^3\sigma$. This lowering of the degree of the interaction can be particularly useful in various contexts, such as the one of the Jacobian Conjecture we are dealing with here.

\section{The Jacobian Conjecture as combinatorial QFT model (the Abdesselam-Rivasseau model)}
\label{subsec-eur}

In \cite{Abdesselam0}, the Jacobian Conjecture was re-expressed as combinatorial QFT model which we call here 
the Abdesselam-Rivasseau model.

Let $n,d\geq 1$, and let 
$F \in \caP_{n,d}$.
The coordinate functions of $F$ can be written as
\begin{equation}
F_i(z)=z_i-\sum_{k=2}^d\sum_{j_1,...,j_k=1}^n
w^{(k)}_{i,j_1...j_k}z_{j_1}...z_{j_k}=: z_i-\sum_{k=2}^d W_i^{(k)}(z)
\ef,
\end{equation}
for $i\leq n$ and $w^{(k)}_{i,j_1...j_k}$ some
coefficients which are nothing but the combinatorial QFT coupling constants (see section \ref{sec:OdimQFT} above).

The partition function
of the Abdesselam-Rivasseau combinatorial QFT model writes
\begin{equation*}
Z(J,K)=\int_{\gC^n}\dd\varphi\dd\varphi^\dag e^{-\varphi^\dag\varphi+\varphi^\dag \sum_{k=2}^d W^{(k)}(\varphi)+J^\dag\varphi+\varphi^\dag K},\nonumber
\end{equation*}
where $J$, $K$ are vectors in $\gC^n$ (and they are nothing but the combinatorial QFT sources, see again section 
\ref{sec:OdimQFT}).
The  
measure of the integral above is 
\beqa
\dd\varphi\dd\varphi^\dag:=\prod_{i=1}^n\frac{\dd {\mathrm{Re}}\varphi_i \, \dd {\mathrm{Im}}\varphi_i}{\pi},
\eeqa
and we have used the standard QFT notations:
\beqa
\varphi^\dag K:=\sum_{i=1}^n \varphi_i^\dag K,  \ \ 
J^\dag\varphi :=\sum_{i=1}^n J^\dag_i 
\varphi_i
\eeqa
Note that the Abdesselam-Rivasseau model has very particular Feynman graphs obtained through perturbative expansion (see the original article \cite{Abdesselam0} for details on this).

Independently of this, setting the coupling constants to zero (the free theory), the partition function is calculated by Gaussian integration:
\begin{equation}
\label{eq-gauss}
\int_{\gC^n}\dd\varphi\dd\varphi^\dag e^{-\varphi^\dag\varphi+J^\dag\varphi+\varphi^\dag K}=e^{J^\dag K}.\nonumber
\end{equation}

Let us now stress the following two crucial points (see again \cite{Abdesselam0} for details):
\begin{enumerate}
\item The partition function $Z$ coincides with the inverse of the Jacobian:
\begin{equation}
Z(0,u)=\det(\partial G(u))= JG(u)=\frac{1}{JF(G(u))}.
\end{equation}
\item The 
inverse $G$ of $F$ corresponds to 
the 
$1-$point 
 correlation 
 function:
\begin{equation}
G_i(u)=\frac{\int_{\gC^n}\dd\varphi\dd\varphi^\dag {\varphi_i} e^{-\varphi^\dag\varphi+\varphi^\dag \sum_{k=2}^d W^{(k)}(\varphi)+\varphi^\dag u}}{\int_{\gC^n}\dd\varphi\dd\varphi^\dag e^{-\varphi^\dag\varphi+\varphi^\dag \sum_{k=2}^d W^{(k)}(\varphi)+\varphi^\dag u}}\label{eq-1point}
\end{equation}
\end{enumerate}

The sets of polynomial functions involved in the Jacobian Conjecture
can be written in this framework:
\begin{align*}
\caJ^{\rm lin}_{n,d}
&=
\{F\in\caP_{n,d}\,|\, Z(0,u)=1,\,\forall u\in\gC^n\},
\\
\caJ_{n,d}
&=
\{F\in\caP_{n,d}\,|\,G_i(u)\text{ given by \eqref{eq-1point} is in
}\caP_{n}\}.
\end{align*}

\section{The intermediate field method for the Abdesselam-Rivasseau model}

As we have seen in 
subsection \ref{sec:inter},
applying the 
combinatorial QFT 
intermediate field method will reduce the
degree $d$ of $F$.
 
We thus {add $n^2$ intermediate fields
$\sigma_{ij}$ ($i,j=1,\ldots, n$) to the model}.
 Indeed, we have, from the general formula
(\ref{eq-gauss}) of Gaussian integration,
\begin{multline}
\label{eq-identfond}
e^{( \varphi^\dag_{i}\varphi_j )
\big(
\sum_{j_2,...,j_d=1}^n
w^{(d)}_{i,j,j_2...j_d}\varphi_{j_2}...\varphi_{j_d} \big) }
\\
=
\int_{\gC^{n^2}}\dd\sigma_{i,j}\dd\sigma^\dag_{i,j} 
e^{-\sigma_{i,j}^\dag\sigma_{i,j}
+\sigma^\dag_{i,j}
\big(
\sum_{j_2,...,j_d=1}^n
w^{(d)}_{i,j,j_2...j_d}\varphi_{j_2}...\varphi_{j_d}
\big)
+(
\varphi^\dag_{i}\varphi_j
)
\sigma_{i,j}}
\end{multline}

We now use the identity \eqref{eq-identfond}, for each pair $(i,j)$,
in the partition function of the model with $n$ dimensions and
degree $d$, in order to re-express the monomials of degree $d$ in
the fields $\varphi$. This leads to:
\begin{multline*}
Z(J,K)=\int_{\gC^n}\dd\varphi\dd\varphi^\dag\int_{\gC^{n^2}}\dd\sigma\dd\sigma^\dag e^{-\varphi^\dag\varphi+\varphi^\dag \sum_{k=2}^{d-1} W^{(k)}(\varphi)+J^\dag\varphi+\varphi^\dag K}\\ e^{\sum_{i,j=1}^n\Big(-\sigma_{i,j}^\dag\sigma_{i,j}+\sigma^\dag_{i,j}\sum_{j_2,...,j_d=1}^n w^{(d)}_{i,j,j_2...j_d}\varphi_{j_2}...\varphi_{j_d}+\varphi^\dag_{i}\varphi_j\sigma_{i,j} \Big)}.
\end{multline*}

\begin{figure}[!tb]
  \centering
\setlength{\unitlength}{50pt}
\begin{picture}(5.1,1.9)
\put(0,0){\includegraphics[scale=1]{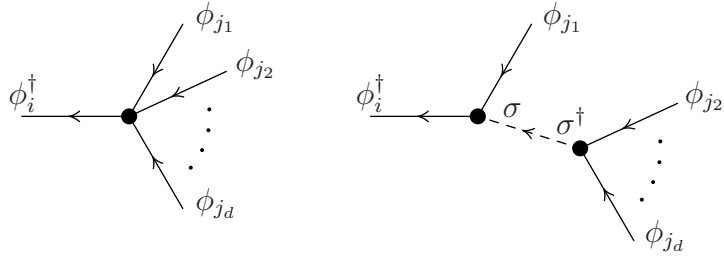}}
\put(0,1.2){$\phi^\dag_i$}
\put(1.4,1.8){$\phi_{j_1}$}
\put(1.7,1.45){$\phi_{j_2}$}
\put(1.4,0.4){$\phi_{j_d}$}
\put(2.6,1.2){$\phi^\dag_i$}
\put(4.0,1.8){$\phi_{j_1}$}
\put(5.063,1.206){$\phi_{j_2}$}
\put(4.763,0.156){$\phi_{j_d}$}
\put(3.7,1.1){$\sigma$}
\put(4.1,0.974){$\sigma^\dag$}
\end{picture}
\caption{Illustration of the intermediate field method. On the left hand side, the interaction term before applying the intermediate field method, and on the right hand side the interaction term after  applying the method.}  
  \label{fig-sigma}
\end{figure}

For a graphical representation of the intermediate field method leading to the model in dimension $n(n+1)$, see Figure~\ref{fig-sigma}.

We define the new vector $\phi$ of $\gC^{n+n^2}$ by
$\phi=(\varphi_1,\ldots,\varphi_n,
\sigma_{1,1}, \ldots, \sigma_{1,n},
\cdots,
\sigma_{n,1}, \ldots, \sigma_{n,n})$.
We further define the interaction coupling constants $\tilde w$ as:
\begin{itemize}
\item for $k=d-1$, we set $\tilde w^{(d-1)}_{i,j,j_2...j_d}:=w^{(d-1)}_{i,j,j_2...j_d}$ and $\tilde w^{(d-1)}_{i\fois n+j,j_2...j_d}=w^{(d)}_{i,j,j_2...j_d}$ with $i,j,j_2,...j_n\leq n$
\item for $k\in\{3,...,d-2\}$, we set $\tilde w^{(k)}_{i,j,j_2...j_k}:=w^{(k)}_{i,j,j_2...j_k}$ with $i,j,j_2,...j_n\leq n$
\item for $k=2$, we set $\tilde w^{(2)}_{i,j,j_2}:=w^{(2)}_{i,j,j_2}$ and $\tilde w^{(2)}_{i,j,i\fois n+j}=1$ with $i,j,j_2\leq n$.
\end{itemize}
The remaining coefficients of $\tilde w$ are set to 0. 

In the same way, the external sources are defined to be
$\tilde J = (J,0)$ and $\tilde K = (K,0)$, where, of course, the
number of extra vanishing coordinates is $n^2$.
It is important to note that these external sources have fewer degrees
of freedom than coordinates ($n$ vs.\ $n(n+1)$).  We also remark that,
for generic $d$, in order to have a relation adapted to induction, it
is crucial to consider an inhomogeneous model, since the intermediate
field method originates terms of degrees $d-1$ and $3$.

The resulting QFT model, after applying the intermediate field method has partition function
\begin{equation*}
Z(J,K)=\int_{\gC^{n+n^2}}\dd\phi\dd\phi^\dag e^{-\phi^\dag\phi+\phi^\dag \sum_{k=2}^{d-1}\tilde W^{(k)}(\phi)+\tilde J^\dag\phi+\phi^\dag \tilde K}
\end{equation*}
and $1-$point correlation function
\begin{equation*}
G_i(u)=\frac{\int_{\gC^{n+n^2}}\dd\phi\dd\phi^\dag \phi_ie^{-\phi^\dag\phi+\phi^\dag \sum_{k=2}^{d-1}\tilde W^{(k)}(\phi)+\phi^\dag \tilde u}}{\int_{\gC^{n+n^2}}\dd\phi\dd\phi^\dag e^{-\phi^\dag\phi+\phi^\dag \sum_{k=2}^{d-1}\tilde W^{(k)}(\phi)+\phi^\dag \tilde u}},
\end{equation*}
for $i\in\{1,\dots,n\}$. 

\bigskip

We have thus showed in a QFT way that the partition function
(resp.\ the one-point correlation function) of the model with dimension
$n\in\gN$ and degree $d\in\gN\setminus\{1,2\}$ is equal to the
partition function (resp.\ the $n$ first coordinates of the one-point
correlation function) of the model with dimension $n(n+1)$ and degree
$d-1$, up to a redefinition of the coupling constant $w\mapsto\tilde
w$ and a trivial redefinition of the external sources. Since, as announced above,  the
partition function corresponds to the inverse of the Jacobian
(resp.\ the one-point correlation function corresponds to the formal
inverse), this gives a QFT proof of
Theorem \ref{lem-main}.
%
This represents 
a reduction
theorem to the quadratic case for Jacobian Conjecture, up to the addition of a
new parameter,
related to the introduction of additional intermediate fields $\sigma$. 
An algebraic proof (not using any QFT arguments) of the Theorem \ref{lem-main} was also derived in the original article \cite{GST}.




\nocite{*}
\bibliographystyle{alpha}
\bibliography{jacobiana}
\label{sec:biblio}

\noindent
{\it\small  
Adrian Tanasa}\\
{\it\small LaBRI,  CNRS UMR 5800, 
universit\'e de Bordeaux}
{\it\small Talence, France, EU}\\
{\it\small H. Hulubei Nat. Inst. Phys. Nucl. Engineering,
Magurele, Romania, EU}\\
{\it\small I. U. F. Paris, France, EU}\\

\end{document}